\documentclass[a4paper,12pt]{amsart}
\usepackage{amsfonts}
\usepackage{amssymb}
\usepackage{ifthen}
\usepackage{graphicx}
\nonstopmode \numberwithin{equation}{section}
\usepackage[usenames]{color}
\usepackage{enumerate}
\usepackage{graphicx}
\usepackage{epstopdf}
\numberwithin{equation}{section}

\theoremstyle{definition}

\newtheorem{remark}[equation]{Remark}

\newtheorem{examples}[equation]{Examples}

\newtheorem{thm}{Theorem}[section]
\newtheorem{lem}[thm]{Lemma}
\newtheorem{cor}[thm]{Corollary}

\newtheorem{cl}{Claim}[section]
\newtheorem{ca}{Case}[section]
\newtheorem{sca}{Subcase}[section]
\newtheorem{scl}[section]{Subclaim}
\newtheorem{conj}[equation]{Conjecture}

\theoremstyle{definition}
\newtheorem{defn}{Definition}[section]
\newtheorem{op}[equation]{Open Problem}
\newtheorem{ques}[equation]{Question}
\newtheorem{exam}[equation]{Example}

\newcounter {own}
\def\theown {\thesection       .\arabic{own}}

\newenvironment{pf}[1][]{%
 \vskip 3mm
 \noindent
 \ifthenelse{\equal{#1}{}}%
  {{\slshape Proof. }}%
  {{\slshape #1.} }%
 }%
{\qed\bigskip}

\newcounter{alphabet}
\newcounter{tmp}

\makeatletter
\newcommand{\reff}[1]{\@ifundefined{r@#1}{}{\setcounter{tmp}{\ref{#1}}\Alph{tmp}}}
\makeatother
\newcommand{\diam}{{\operatorname{diam}}}

\newcommand{\dist}{{\operatorname{dist}}}

\def\be{\begin{equation}}
\def\ee{\end{equation}}

\newcommand{\bee}{\begin{enumerate}}
\newcommand{\eee}{\end{enumerate}}

\newcommand{\blem}{\begin{lem}}
\newcommand{\elem}{\end{lem}}
\newcommand{\bthm}{\begin{thm}}
\newcommand{\ethm}{\end{thm}}
\newcommand{\bcor}{\begin{cor}}
\newcommand{\ecor}{\end{cor}}
\newcommand{\beg}{\begin{exam}}
\newcommand{\eeg}{\end{exam}}
\newcommand{\begs}{\begin{examples}}
\newcommand{\eegs}{\end{examples}}
\newcommand{\bdefe}{\begin{defn}}
\newcommand{\edefe}{\end{defn}}
\newcommand{\bprob}{\begin{prob}}
\newcommand{\eprob}{\end{prob}}
\newcommand{\bques}{\begin{ques}}
\newcommand{\eques}{\end{ques}}
\newcommand{\bei}{\begin{itemize}}
\newcommand{\eei}{\end{itemize}}
\newcommand{\bcon}{\begin{conj}}
\newcommand{\econ}{\end{conj}}
\newcommand{\bop}{\begin{op}}
\newcommand{\eop}{\end{op}}

\newcommand{\bca}{\begin{ca}}
\newcommand{\eca}{\end{ca}}
\newcommand{\bsca}{\begin{sca}}
\newcommand{\esca}{\end{sca}}

\newcommand{\bcl}{\begin{cl}}
\newcommand{\ecl}{\end{cl}}

\newcommand{\bscl}{\begin{scl}}
\newcommand{\escl}{\end{scl}}

\newcommand{\bcons}{\begin{conjs}}
\newcommand{\econs}{\end{conjs}}
\newcommand{\bprop}{\begin{propo}}
\newcommand{\eprop}{\end{propo}}
\newcommand{\er}{\end{rem}}
\newcommand{\brs}{\begin{rems}}
\newcommand{\ers}{\end{rems}}
\newcommand{\bo}{\begin{obser}}
\newcommand{\eo}{\end{obser}}
\newcommand{\bos}{\begin{obsers}}
\newcommand{\eos}{\end{obsers}}
\newcommand{\bpf}{\begin{pf}}
\newcommand{\epf}{\end{pf}}
\newcommand{\ba}{\begin{array}}
\newcommand{\ea}{\end{array}}
\newcommand{\beq}{\begin{eqnarray}}
\newcommand{\beqq}{\begin{eqnarray*}}
\newcommand{\eeq}{\end{eqnarray}}
\newcommand{\eeqq}{\end{eqnarray*}}

\newcounter{minutes}
\divide\time by 60
\newcounter{hours}
\multiply\time by 60 \addtocounter{minutes}{-\time}
\begin{document}

\bibliographystyle{amsplain}

\title[The weak min-max property in Banach spaces]{The weak min-max property in Banach spaces}

%%%%%%%% BEGIN TIMESTAMP
%\def\thefootnote{}
%\footnotetext{ \texttt{\tiny File:~\jobname .tex,
%          printed: \number\year-\number\month-\number\day,
%          \thehours.\ifnum\theminutes<10{0}\fi\theminutes}
%} \makeatletter\def\thefootnote{\@arabic\c@footnote}\makeatother
%%%%%%%% END TIMESTAMP

\author[Z. Ouyang]{Zhengyong Ouyang}
\address{Zhengyong Ouyang, School of Mathematics and Big Data, Foshan University,  Foshan, Guangdong 528000, People's Republic
of China} \email{oyzy1128@126.com}

\author{Antti Rasila}
\address{Antti Rasila,
Mathematics with Computer Science Program,
Guangdong Technion -- Israel Institute of Technology
241 Daxue Road, Shantou 515063, Guangdong, People's Republic of China
and
Department of Mathematics,
Technion -- Israel Institute of Technology
Haifa 32000, Israel} \email{antti.rasila@gtiit.edu.cn; antti.rasila@iki.fi}

\author{Tiantian Guan$^{~\mathbf{*}}$}
\address{Tiantian Guan, School of Mathematics and Big Data, Foshan University,  Foshan, Guangdong 528000, People's Republic
of China} \email{ttguan93@163.com}

\begin{abstract}
In this paper, we investigate the relationship between the weak min-max property and the diameter uniformity of domains in Banach spaces with dimensions at least $2$. As an application, we show that diameter uniform domains are invariant under relatively quasim\"obius mappings.
\end{abstract}

\date{}
 \keywords{Uniform domain, weak min-max property,
(relatively) quasim\"{o}bius mapping.\\
${}^{\mathbf{*}}$ Corresponding author}

%\thanks{
%Qingshan Zhou was supported by NNSF of China (Nos. 11901090), and by Department of Education of Guangdong Province, China (Grant Nos. 2018KQNCX285 and 2018KTSCX245).
\thanks{Antti Rasila was supported by NNSF of China (No. 11971124), and NSF of the Guangdong Province (No. 2021A151010326). Tiantian Guan was supported by the Guangdong Basic and Applied Basic Research Foundation under Grant (No. 2021A1515012289), and Research Fund of Guangdong-Hong Kong-Macao Joint Laboratory for Intelligent Micro-Nano Optoelectronic Technology (No. 2020B1212030010).}

%\maketitle{} \pagestyle{myheadings} \markboth{}{Quasim\"obius invariance of uniform domains}

\maketitle

\section{Introduction and main results}\label{sec-1}

 Uniform domains in $\mathbb{R}^n$ were introduced by Martio and Sarvas in 1978 \cite{MS}, who used them in studying injectivity theorems in $\mathbb{R}^n$. Since then, uniform domains have played an important role in the classical function theory, quasiconformal mappings, and many other fields of modern mathematical analysis (see e.g. \cite{GP,Jo80,Vai-0,Vai,Vai-5}). Several different characterizations of uniform domains have been established, see for example \cite{GH, GO, LPZ, LVZ2, LVZ3, Martio-80}.

In particular, Martio  \cite{Martio-80} showed that both the class of diameter uniform domains and the class of $\delta$-uniform domains are equivalent to the class of uniform domains in $\mathbb{R}^n$. Gehring and Hag \cite{GH} introduced the min-max property of a domain $G\subsetneq \mathbb{R}^n$ involving the properties of hyperbolic geodesic in $\mathbb{B}^n$. By using this property, they showed other characteristics of uniform domains in $\mathbb{R}^n$. Note that the above classes of domains in $\mathbb{R}^n$ are equivalent. The definitions of these classes of domains are given in Section \ref{sec-2}.

Through this paper, let $E$ and $E'$ be real Banach spaces with dimension at least $2$. For domains in the setting of Banach spaces, V\"ais\"al\"a \cite{Vai04} showed that the class of distance uniform domains is equivalent to the diameter uniform domains. Furthermore, it is known that uniform domains are diameter uniform. For the converse, Rasila and Zhou \cite{ZR} showed that a diameter uniform domain is uniform if the domain satisfies a natural condition, which is as follows:

\begin{thm}$($\cite[Theorem 1.2]{ZR}$)$\label{ZR-thm1.2}
Suppose that $G\subsetneq E$ is a domain. Then $G$ is $A$-uniform if and only if $G$ is diameter $A_1$-uniform and $\psi$-natural, where $A$, $A_1$ and $\psi$ depend on each other.
\end{thm}

Concerning on uniform domains in Banach spaces $E$, Rasila and Zhou studied the following implications:

\begin{thm}$($\cite[Theorem 1.1]{ZR}$)$\label{ZR-thm1.1}
 Suppose that $G\subsetneq E$ is a domain. If $G$ is $A$-uniform, then $G$ satisfies the $c$-min-max property, $G$ is a diameter $A_1$-uniform domain, where $c=c(A)$, $A_1=A_1(A)$, and $c=c(A)$ means that $c$ depends only on $A$. Furthermore, $G$ is a diameter $A_1$-uniform domain if and only if $G$ is $\delta$-uniform for some $0<\delta<1$, where $A_1$ and $\delta$ depend on each other.
\end{thm}

It is reasonable to consider the question whether diameter uniform domains satisfy the min-max property. To deal with this question, we introduce a condition called weak min-max property, and use it to characterize diameter domains in Banach spaces.

\bdefe Let $G \subsetneq E$ be a domain. We say that $G$ has the {\it weak min-max property} if there exists a family of curves $\Gamma$ in $G$ and a constant $c\geq 1$ such that any $x_1$ and $x_2$ in $G$ can be joined by a curve $\gamma\in \Gamma$ satisfying
\be%\label{r-1} 
c^{-1}\min_{j=1,2} |x_j-y| \leq |x-y| \leq c\max_{j=1,2} |x_j-y|, \ee
for all $x\in \gamma$ and all $y\in \partial G$.
\edefe

\begin{thm}\label{thm-1}
Suppose that $G\subsetneq E$ is a domain. Then $G$ has the weak $c$-min-max property if and only if $G$ is diameter $A$-uniform, where $c$ and $A$ depend only on each other.
\end{thm}

It follows from \cite[Theorem 6.26]{Vai-2} that uniform domains are preserved under quasim\"obius mappings. In \cite[Theorem 1.3]{ZR}, Rasila and Zhou recently proved that relatively quasim\"obius mappings map uniform domains to diameter uniform domains. As an application of Theorem \ref{thm-1}, we study the relative quasim\"obius invariance of diameter uniform domains.

\begin{thm}\label{thm-2}
Suppose that $G\subsetneq E$ and $G'\subsetneq E'$ are domains. Suppose that $G$ is diameter $A$-uniform, and that a homeomorphism $f:\, \overline{G}\to \overline{G'}$ is $\theta$-quasim\"obius relative to $\partial G$ and maps $G$ onto $G'$. Then $G'$ is diameter $A_1$-uniform, where $A_1=A_1(A, \theta)$.
\end{thm}

Using Theorems \ref{thm-1} and \ref{thm-2}, we see that the weak min-max property is invariant under relatively quasim\"obius mappings.

\begin{cor}\label{cor-0}
Suppose that $G\subsetneq E$ and $G'\subsetneq E'$ are domains. Assume that $G$ has the weak $c$-min-max property, and that  a homeomorphism $f:\, \overline{G}\to \overline{G'}$ is $\theta$-quasim\"obius relative to $\partial G$ and maps $G$ onto $G'$. Then $G'$ satisfies the weak $c_1$-min-max property, where $c_1=c_1(c, \theta)$.
\end{cor}

In \cite[Theorem 1.3]{ZR}, it was shown that relatively quasim\"obius mappings map uniform domains onto uniform domains if the target domains are natural. By Corollary \ref{cor-0} and Theorem \ref{ZR-thm1.2}, we have the following corollary.

\begin{cor}\label{cor-1}
Suppose that $G\subsetneq E$ and $G'\subsetneq E'$ are domains. Suppose that $G$ is diameter $A$-uniform, $G'$ is $\psi$-natural, and that $f:\, G\to G'$ is $\theta$-quasim\"obius relative to $\partial G$. Then $G$ is $A_1$-uniform, where $A_1=A_1(A, \theta, \psi)$.
\end{cor}

%\br The min-max property for domains in $\mathbb{R}^n$ was introduced by Gehring and Hag in \cite{GH}. They extended the properties of hyperbolic geodesic in $\mathbb{B}^n$ to more general domains, and studied the relationship between this property, the uniformity and quasiconformal extension property. In \cite{Martio-80}, Martio introduced the concept of the $\delta$-uniformity property in terms of the cross-ratio of four points for domains in $\mathbb{R}^n$. By using this condition, he obtained certain general properties of uniform domains.
%\er

%\br The definition of $\psi$-natural domains is given in Definition \ref{n-1}.  We note that every proper domain in $\mathbb{R}^n$ is $\psi$-natural with $\psi=\psi(n)$ (see \cite[Corollary 2.18]{Vu1}). In an infinite dimensional Hilbert space, the broken tube construction  in \cite{Vai04} provides an example of a domain, which is not natural. Moreover, it is diameter $c$-uniform but not $c_1$-uniform for any constant $c_1\geq 1$.
%\er

The rest of this paper is organized as follows. In Section \ref{sec-2}, we recall necessary definitions and preliminary results. The proof of Theorem \ref{thm-1} is given in Section \ref{sec-3}. Section \ref{sec-4} is devoted to the proofs of Theorem \ref{thm-2} and Corollary \ref{cor-0}.

\section{Preliminaries}\label{sec-2}

\subsection{\bf Notation.}
Following the notation and terminology of \cite{HLVW,Vai-5}, the norm of a vector $x$ in $E$ is written as $|x|$, and for every pair of points $x$, $y$ in $E$, the distance between them is denoted by $|x-y|$, the closed line segment with endpoints $x$ and $y$ by $[x,y]$. %The one-point extension of $E$ is the Hausdorff space $\dot{E}=E\cup\{\infty\}$, where the neighborhoods of $\infty$ are the complements of closed bounded sets of $E$.

For a set $A$ in $E$, we use $\overline{A}$ to denote the completion of $A$ and $\partial A=\overline{A}\setminus A$ to be its norm boundary. For a bounded set $A$ in $E$, $\diam (A)$ is the diameter of $A$. %Let
%$$B(x,r)=\{ z\in E\,|\; |z-x|<r\},\; \overline{B}(x,r)=\{ z\in E\,|\; |z-x|\leq r\},$$
%and $S(x,r)=\{ z\in E\,|\; |z-x|=r\}$.
 A curve is a continuous function $\gamma:$ $[a,b]\to E$. The length of $\gamma$ is defined by
$$\ell(\gamma)=\sup\Big\{\sum_{i=1}^{n}|\gamma(t_i)-\gamma(t_{i-1})|\Big\},$$
where the supremum is taken over all partitions $a=t_0<t_1<t_2<\ldots<t_n=b$.

\subsection{\bf  Uniform domains}\label{sub-domain}
\bdefe\label{d-du} Let $G \subsetneq E$ be a domain and let $A\geq 1$. Then $G$ is called {\it diameter $A$-uniform}, if each pair of points $x_{1},x_{2}$ in $G$ can be joined by a curve $\alpha$ in $G$ satisfying:
\begin{enumerate}
\item\label{en-1} $\min\limits_{j=1,2}\;\diam\; (\alpha [x_j, x])\leq A\,d_{G}(x)
$ for all $x\in \alpha$, and
\item\label{en-2} $\diam \,(\alpha) \leq A\,|x_{1}-x_{2}|$,
\end{enumerate}
where $d_G(x)=\dist(x, \partial G)$. Moreover, the curve $\alpha$ satisfying the above conditions is said to be a {\it diameter uniform curve}.
\edefe
The domain $G\subsetneq E$ is called {\it $A$-uniform} if the diameter conditions in Definition \ref{d-du} are replaced by the length conditions. We say that $G$ is {\it distance $A$-uniform} if the first condition in Definition \ref{d-du} is replaced by
$$\min\limits_{j=1,2}|x_j-x|\leq A d_G(x),$$
for all $x\in \alpha$.

\begin{remark} It is not difficult to see from the definitions that an $A$-uniform domain is diameter $A$-uniform and a diameter $A$-uniform domain is distance $A$-uniform. For the converse, diameter uniform domains are equivalent to distance uniform domain in Banach spaces \cite[Page 10]{Vai04}. It follows from \cite[Theorem 4.5]{Martio-80} that a diameter $A$-uniform domain in $\mathbb{R}^n$ is $A_1$-uniform, where $A_1$ depends only on $A$ and $n$. Note that in infinite dimensional Hilbert spaces, V\"ais\"al\"a \cite{Vai04} constructed a so-called broken tube domain, which is diameter $A$-uniform but not $A_1$-uniform for any constant $A_1\geq 1$.
\end{remark}

\bdefe Let $G \subsetneq E$ be a domain and let $0<\delta<1$. Then $G$ is called $\delta$-{\it uniform} if each pair of points $x_{1},x_{2}$ in $G$ can be joined by a curve $\alpha$ in $G$ such that the cross ratio
\be\label{r-2} \tau(x,x_i,y,x_j)=\frac{|x-y|}{|x-x_i|} \cdot \frac{|x_i-x_j|}{|x_j-y|}\geq \delta,\;\;\;i\neq j\in\{1,2\}, \ee
for all $x\in \alpha\setminus\{x_1,x_2\}$ and $y\in E\setminus G$.
\edefe

\bdefe\label{n-1}  Let $\psi:[0,\infty)\to [0,\infty)$ be an increasing function. A domain $G \subsetneq E$ is called $\psi$-{\it natural} if
$$k_G(A)\leq \psi(r_G(A))$$
for every nonempty connected set $A\subseteq G$ with $r_G(A)<\infty$, where $$r_G(A)=\sup\Big\{\frac{|x-y|}{\min\{d_G(x), d_G(y)\}}\,|\, x\in A, y\in A\Big\},$$
 and $k_G(A)$ is the diameter of $A$ in the quasihyperbolic metric.
\edefe
\begin{remark}
The class of natural domains is fairly large, for example all proper domains in $\mathbb{R}^n$ is $\psi$-natural with $\psi=\psi(n)$, see \cite[Corollary 2.18]{Vu1}. In an infinite dimensional Hilbert space the broken tube construction given in \cite{Vai04} provides an example of a domain, which is not natural.
\end{remark}

\subsection{Min-max property}

\bdefe Let $G \subsetneq E$ be a domain. We say that $G$ has the {\it min-max property} if there exists a family of curves $\Gamma$ in $G$ and a constant $c\geq 1$ such that any pair of points in $G$ can be joined by a curve $\gamma\in \Gamma$ with the property that
\be\label{r-1} c^{-1}\min_{j=1,2} |x_j-y| \leq |x-y| \leq c\max_{j=1,2} |x_j-y|, \ee
for each ordered triple of points $x_1,x,x_2\in \gamma$ and each $y\in \partial G$.
\edefe

\subsection{\bf Quasisymmetric mappings and quasim\"obius mappings}\label{sub-1}

By a triple in $X$ we mean an ordered sequence $T=(x,y,z)$ of three
distinct points in $X$. The ratio of $T$ is the number
$$\rho(T)=\frac{|y-x|}{|z-x|}.$$
If $f: X\to Y$ is  an injective
map, then the image of a triple  $T=(x,y,z)$  is the triple
$f(T)=(f(x),f(y),f(z))$.

Suppose that  $A\subseteq X$. A triple  $T=(x,y,z)$ in $X$ is said to
be a triple in the pair $(X, A)$ if $x\in A$ or if $\{y,z\}\subseteq
A$. Equivalently, both $|y-x|$ and $|z-x|$ are distances from a
point in $A$.

\bdefe \label{def1-0} Let $X$ and $Y$ be two metric spaces, and let
$\eta: [0, \infty)\to [0, \infty)$ be a homeomorphism. Suppose
$A\subset X$. A homeomorphism $f: X\to Y$ is said to be {\it
$\eta$-quasisymmetric} relative to $A$, if $\rho(f(T))\leq \eta(\rho(T))$  for each triple $T$ in
$(X,A)$. \edefe

Note that a homeomorphism $f: X\to Y$ is $\eta$-quasisymmetric relative to  $A$ if
and only if $\rho(T)\leq t$ implies that $\rho(f(T))\leq \eta(t)$
for each triple $T$ in $(X,A)$ and $t\geq 0$ (cf. \cite{Vai-5}).
Obviously, a quasisymmetric mapping relative to  $X$ is equivalent to usual quasisymmetry.

%A quadruple in $X$ is an ordered sequence $Q=(x,y,z,w)$ of four
%distinct points in $X$. The cross ratio of $Q$ is defined to be the
%number
%$$\tau(Q)=\tau(x,y,z,w)=\frac{|x-z|}{|x-y|}\cdot\frac{|y-w|}{|z-w|}.$$
Observe that the definition of the cross ratio can be extended in
the usual manner to the case where one of the points is
$\infty$. For example,
$$\tau(x,y,z,\infty)= \frac{|x-z|}{|x-y|}.$$
Let $X$ be a metric space and $\dot{X}=X\cup \{\infty\}$. If $X_0 \subseteq \dot{X}$ and $f: X_0\to \dot{Y}$
is an injective map, then the image of a quadruple $Q=(x, y, z, w)$ in $X_0$ is the
quadruple $f(Q)=(f(x),f(y),f(z),f(w))$. Suppose that $A\subseteq X_0$. We say
that a quadruple $Q=(x,y,z,w)$ in $X_0$ is a quadruple in the pair
$(X_0, A)$ if $\{x,w\}\subseteq A$ or $\{y,z\}\subseteq A$.
Equivalently, all four distances in the definition of $\tau(Q)$ are
 distances from a point in $A$.

\bdefe \label{def2'}
Let ${X}$ and ${Y}$ be two metric
spaces, $X_1\subset \dot{X}$, and $Y_1\subset \dot{Y}$, and let $\eta: [0, \infty)\to [0, \infty)$ be a
homeomorphism. Suppose $A\subseteq \dot{X}$. A homeomorphism $f:
X_1\to Y_1$ is said to be {\it $\eta$-quasim\"obius}
relative to $A$, if the inequality
$\tau(f(Q))\leq \eta(\tau(Q))$ holds for each quadruple in $(X,A)$.
\edefe

Apparently, $\eta$-quasim\"obius relative to  $X$ is equivalent to $\eta$-quasim\"obius.

 %%%%%%%%%%%%%%%%%%%%%%%%%%%%%%%%%%%%%%%%%
 %%%%%%%%%%%%%%%%%%%%%%%%%%%%%%%%%%%%%%%%%
\section{Proof of Theorem \ref{thm-1}}\label{sec-3}
 %%%%%%%%%%%%%%%%%%%%%%%%%%%%%%%%%%%%%%%%%
 %%%%%%%%%%%%%%%%%%%%%%%%%%%%%%%%%%%%%%%%%
Assume that $G\subsetneq E$ is a domain.
First, we show the sufficiency part and with the assumption that $G$ is diameter $A$-uniform.

Fix $x_1, x_2\in G$. By the assumption, there is a diameter uniform curve $\alpha$ joining $x_1$ and $x_2$ in $G$. For all $x\in \gamma$ and for all $y\in \partial G$, we have
$$\min\limits_{j=1,2}|x_j-x|\leq \min\limits_{j=1,2}\diam(\alpha[x, x_j])\leq Ad_G(x)\leq A|x-y|,$$
which implies
\be\label{e3-1}
\min\limits_{j=1,2}|x_j-y|\leq \min\limits_{j=1,2}|x_j-x|+|x-y|\leq (A+1)|x-y|.
\ee
Moreover, we see from the diameter uniformity of $\alpha$ that
\beqq
|x-y|&\leq& \diam(\alpha)+|x_1-y|\\
&\leq& A|x_1-x_2|+|x_1-y|\\
&\leq& (2A+1)\max\limits_{j=1,2}|x_j-y|.
\eeqq
This, together with \eqref{e3-1}, shows that
$$c^{-1}\min\limits_{j=1,2}|x_j-y|\leq |x-y|\leq c\max\limits_{j=1,2}|x_j-y|,$$
where $c=2A+1$.

Now, we prove the necessity. Assume that $G$ has the weak $c$-min-max property. Fix $x_1, x_2\in G$. Without loss of generality, we may assume that
\be\label{e3-2}
|x_1-x_2|>\frac{1}{2} d_G(x_1).
\ee
 Indeed, if $|x_1-x_2|\leq 1/2 d_G(x_1),$ then it is easy to see that the line segment $[x_1, x_2]$ is the desired diameter uniform curve.

 Join $x_1$ and $x_2$ by a curve $\gamma$ satisfying the weak min-max property, and fix $y\in \partial G$ with $|x_1-y|\leq 2d_G(x_1)$. Then it follows from \eqref{e3-2} that for all $x\in \gamma$, we have
 \beqq\label{e3-3}
 |x-y|\leq c\max\limits_{j=1,2}|x_j-y|\leq c(|x_1-x_2|+|x_1-y|)\leq 5c|x_1-x_2|,
 \eeqq
which yields
\be\label{e3-4}
\diam(\gamma)\leq 10c|x_1-x_2|.
\ee
For $x\in \gamma$, choose $z\in \partial G$ with $|x-z|\leq 2d_G(x)$. Thus by the weak min-max property, we have
\be\label{e3-6}
\min\limits_{j=1,2}|x_j-x|\leq \min\limits_{j=1,2}|x_j-z|+|z-x|\leq (c+1)|z-x|\leq 2(c+1)d_G(x).
\ee
By \eqref{e3-4} and \eqref{e3-6}, we obtain that $G$ is distance $A$-uniform with $A=10c$. Therefore, it follows from \cite[Page 10]{Vai04} that $G$ is diameter $A_1$-uniform, where $A_1=A_1(c)$.
%where $A_1$ depends only on $c$.
\qed

 %%%%%%%%%%%%%%%%%%%%%%%%%%%%%%%%%%%%%%%%%
 %%%%%%%%%%%%%%%%%%%%%%%%%%%%%%%%%%%%%%%%%
\section{Proofs of Theorem \ref{thm-2} and Corollary \ref{cor-1}}\label{sec-4}
 %%%%%%%%%%%%%%%%%%%%%%%%%%%%%%%%%%%%%%%%%
 %%%%%%%%%%%%%%%%%%%%%%%%%%%%%%%%%%%%%%%%%
 Before giving the proof of Theorem \ref{thm-2}, we show that the weak min-max property is invariant under relatively quasisymmetric mappings.

\begin{lem}\label{lem4-1}
Let $G\subsetneq E$ and $G'\subsetneq E'$ be domains. Suppose that $G$ has the weak $c$-min-max property, and that a homeomorphism $f:\,\overline{G}\to \overline{G'}$ is $\eta$-quasisymmetric relative to $\partial G$ and maps $G$ onto $G'$, then $G'$ satisfies the weak $c_1$-min-max property, where $c_1=c_1(c, \eta)$.
\end{lem}
\bpf
Fix $x_1', x_2'\in G'$. Then we find two points $x_1,x_2\in G$ such that $x_j'=f(x_j)$ for $j=1, 2$. By the min-max property of $G$, there is a curve $\gamma$ in $G$ joining $x_1$ and $x_2$ such that for all $x\in \gamma$ and for all $y\in \partial G$,
\be\label{eq4-1}
c^{-1}\min\limits_{j=1,2}|x_j-y|\leq |x-y|\leq c\max\limits_{j=1,2}|x_j-y|
\ee
Without loss of generality, we may assume that $|x_1-y|\leq |x_2-y|$. Hence, by \eqref{eq4-1}, and by the relative quasisymmetry of $f$, we have
\beqq\label{eq4-2}
|x_1'-f(y)|=|f(x_1)-f(y)|\leq \eta(c)|f(x)-f(y)|
\eeqq
and
\beqq
|f(x)-f(y)|\leq \eta(c)|f(x_2)-f(y)|=\eta(c)|x_2'-f(y)|,
\eeqq
which imply
$$\eta(c)^{-1}\min\limits_{j=1,2}|x_j'-f(y)|\leq |f(x)-f(y)|\leq \eta(c)\max\limits_{j=1,2}|x_j'-f(y)|.$$
It follows that $G'$ has the weak $c_1$-min-max property with $c_1=\eta(c)$.
\epf

V\"ais\"al\"a showed that the inversions are quasim\"obius and the relationship between quasisymmetric mappings and quasim\"obius mappings, which are as follows:

\begin{lem}(\cite[Page 220]{Vai-0})\label{Vai-0P220}
Let $u(x)$ be the inversion in the unit sphere $\mathbb{S}$, where
$$u(x)=\frac{x}{|x|^2}$$
for $x\in E\setminus \{o\}$ and $o$ is the origin. Then $u^{-1}=u$, and $u$ is $\theta_0$-quasim\"obius with the control function $\theta_0(t)=81t$ for $t\in [0, \infty)$.
\end{lem}

\begin{thm}(\cite[Theorem 3.10]{Vai-0})\label{Vai-0thm310}
Let $X$ and $Y$ be two metric spaces. Suppose that $X$ is unbounded and that $f:X\to Y$ is $\theta$-quasim\"obius. Then $f$ is quasisymmetric if and only if $f(x)\to \infty$ as $x\to \infty$. In this case, $f$ is $\theta$-quasisymmetric.
\end{thm}

Now by Theorem \ref{thm-1} and Lemma \ref{lem4-1}, we know that diameter uniform domains are preserved under relatively quasisymmetric mappings.

\begin{cor}\label{cor-2}
Let $G\subsetneq E$ and $G'\subsetneq E'$ be domains. Suppose that $G$ is diameter $A$-uniform, and that  $f:\,\overline{G}\to \overline{G'}$ is $\eta$-quasisymmetric relative to $\partial G$ homeomorphism maps $G$ onto $G'$. Then $G'$ is diameter $A_2$-uniform, where $A_2=A_2(A, \eta)$.
\end{cor}

Let $u(x)$ be the inversion in the unit sphere $\mathbb{S}$ as defined in Lemma \ref{Vai-0P220}. It follows from Lemma \ref{Vai-0P220} that $u$ and $u^{-1}$ are $\theta_0$-quasim\"obius with the control function $\theta_0(t)=81t$, $t\in [0,\infty)$.

Next, we show that diameter uniformity is invariant under an inversion transformation.

\begin{lem}\label{lem4-2}
Suppose that $G$ is a diameter $A$-uniform domain of $E$. Then $u(G)$ is diameter $A_2$-uniform with $A_3=A_3(A)$.
\end{lem}
\bpf
By Theorem \ref{ZR-thm1.1}, the diameter uniformity of $G$ gives that $G$ is $\delta$-uniform with $\delta=\delta(A)$. Because $u$ is $\theta_0$-quasim\"obius with the control function $\theta_0(t)=81t$, $t\in [0,\infty)$, a simple computation shows that $u(G)$ is $\delta_1$-uniform with  $\delta_1=\delta_1(A)$. Using Theorem \ref{ZR-thm1.1}, we know that $u(G)$ is diameter $A_3$-uniform with $A_3=A_3(A)$.
\epf

\subsection*{Proof of Theorem \ref{thm-2}}
By composing translations if necessary, we may assume that $o\in \partial G$ and $o\in \partial G'$. Denote
$$g:=\,u\circ f\circ u^{-1}:\,u(G)\to u(G').$$
It follows from Lemma \ref{Vai-0P220} that $u$ and $u^{-1}$ are $\theta_0$-quasim\"obius with the control function $\theta_0(t)=81t$, $t\in [0,\infty)$. Hence, we obtain that $g$ is $\theta_1$-quasim\"obius relative to $ u(\partial G)$ with $g(x)\to \infty$ as $x\to \infty$, where the control function $\theta_1(t)=81\theta(81t)$, $t\in [0,\infty)$. By Theorem \ref{Vai-0thm310}, $g$ is $\theta_1$-quasisymmetric relative to $u(\partial G)$ with $\theta_1=\theta_1(\theta)$.

On the other hand, by the diameter uniformity of $G$, Lemma \ref{lem4-2}, yields that $u(G)$ is diameter $A_3$-uniform. Therefore, we see form Corollary \ref{cor-2} that $u(G')$ is diameter $A_4$-uniform with $A_4=A_4(A, \theta)$. Because $u^{-1}=u$, Lemma \ref{lem4-2} yields that $G'$ is diameter $A_1$-uniform with $A_1=A_1(A, \theta)$.
\qed

\subsection*{Proof of Corollary \ref{cor-1}.} By Theorem \ref{thm-2}, we obtain that $G'$ is diameter $A'$-uniform, where $A'=A'(A, \theta)$. It follows from Theorem \ref{ZR-thm1.2} that $G'$ is $A_2$-uniform, where $A_2=A_2(A, \theta, \psi)$.
\qed

%{\bf Acknowledgement.}
%%%%%%%%%%%%%%%%%%%
%%%%%%%%%%%%%%%%%%%

%{\bf Acknowledgement.} The authors are indebted to the referees for the valuable suggestions.


\begin{thebibliography}{99}

\bibitem{GH}  {\sc F. W. Gehring and K. Hag}, Remarks on uniform and quasiconformal extension domains, \textit{Complex Variables Theory Appl.} {\bf 9} (1987), 175--188.

\bibitem{GO}  {\sc F. W. Gehring and B. G. Osgood}, Uniform domains and the quasi-hyperbolic metric, \textit{J. Anal. Math.} {\bf 36} (1979), 50--74.

\bibitem{GP}  {\sc F. W. Gehring and B. P. Palka}, Quasiconformally homogeneous domains, \textit{J. Anal. Math.} {\bf 30} (1976), 172--199.

%\bibitem{HL}  {\sc X. Huang and J. Liu}, Quasihyperbolic metric and quasisymmetric mappings in metric spaces, {\it Trans. Amer. Math. Soc.} {\bf 367} (2015), 6225--6246.

\bibitem{HLVW}  {\sc M. Huang, Y. Li, M. Vuorinen, and X. Wang}, On quasim\"obius maps in real Banach spaces, \textit{Israel J. Math.} {\bf 198} (2013), 467--486.

%\bibitem{HRWZ}  {\sc M. Huang, A. Rasila, X. Wang, and Q. Zhou}, Semisolidity and locally weak quasisymmetry of homeomorphisms in metric spaces, \textit{ Stud. Math.} {\bf 242} (2018), 267--301.

%\bibitem{Jo61} {\sc F. John}, Rotation and strain,  \textit{Comm. Pure. Appl. Math.} {\bf 14} (1961), 391--413.

\bibitem{Jo80} {\sc P. W. Jones}, Extension theorems for BMO, \textit{Indiana Univ. Math. J.}  {\bf 29} (1980), 41--66.

\bibitem{LPZ} {\sc Y. Li, S. Ponnusamy, and Q. Zhou}, Sphericalization and flattening preserve uniform domains in non-locally compact metric spaces, \textit{J. Aust. Math. Soc.} DOI:10.1017/S1446788719000582  (2020),  1--22.

\bibitem{LVZ2} {\sc Y. Li, M. Vuorinen, and Q. Zhou}, Weakly quasisymmetric maps and uniform spaces, \textit{Comput. Methods Funct. Theory}  {\bf 18} (2018), 689--715.

%\bibitem{LVZ17} {\sc Y. Li, M. Vuorinen, and Q. Zhou}, Characterizations of John spaces, \textit{Monatsh. Math.}  {\bf 188} (2019), 547--559.

\bibitem{LVZ3} {\sc Y. Li, M. Vuorinen, and Q. Zhou}, Apollonian metric, uniformity and Gromov hyperbolicity, \textit{Complex Var. Elliptic Equ.}  {\bf 65} (2020), 215--228.

\bibitem{Martio-80}  {\sc O. Martio}, Definitions of uniform domains, \textit{Ann. Acad. Sci. Fenn. Ser. A I Math.} {\bf 5} (1980), 197--205.

\bibitem{MS}  {\sc O. Martio and J. Sarvas}, Injectivity theorems in plane and space, \textit{Ann. Acad. Sci. Fenn. Ser. A I Math.} {\bf 4} (1978), 383--401.

%\bibitem{rt2} {\sc A. Rasila and J. Talponen,} On Quasihyperbolic Geodesics in Banach Spaces,  \textit{Ann. Acad. Sci. Fenn. Ser. A I Math.} {\bf 39} (2014), 163--173.

%\bibitem{TV}  {\sc P. Tukia and J. V\"{a}is\"{a}l\"{a}}, Quasisymmetric embeddings of metric spaces,
%\textit{Ann. Acad. Sci. Fenn. Ser. A I Math.,} {\bf 5} (1980), 97--114.

\bibitem{Vai-0}  {\sc J. V\"{a}is\"{a}l\"{a}}, Quasim\"obius maps, \textit{J. Anal. Math.} {\bf 44} (1984/85), 218--234.

\bibitem{Vai}  {\sc J. V\"{a}is\"{a}l\"{a}}, Uniform domains, \textit{Tohoku Math. J.} {\bf 40} (1988), 101--118.

%\bibitem{Vai-1}  {\sc J. V\"{a}is\"{a}l\"{a}}, Free quasiconformality in Banach spaces. I, \textit{Ann. Acad. Sci. Fenn. Ser. A I Math.} {\bf 15} (1990), 355--379.

\bibitem{Vai-2}  {\sc J. V\"{a}is\"{a}l\"{a}}, Free quasiconformality
in Banach spaces. II, \textit{Ann. Acad. Sci. Fenn. Ser. A I Math.}
{\bf 16} (1991), 255-310.

\bibitem{Vai-3}  {\sc J. V\"{a}is\"{a}l\"{a}}, Free quasiconformality
in Banach spaces. III, \textit{Ann. Acad. Sci. Fenn. Ser. A I Math.}
{\bf 17} (1992), 393--408.

%\bibitem{Vai-4}  {\sc J. V\"{a}is\"{a}l\"{a}}, Free quasiconformality in Banach spaces. IV, \textit{Analysis and Topology, 697--717, World Sci. Publ., River Edge, N. J.,} 1998.

%\bibitem{Vai4}  {\sc J. V\"{a}is\"{a}l\"{a}}, Relatively and inner
% uniform domains, \textit{Conformal Geom. Dyn.,} {\bf 2} (1998), 56--88.


\bibitem{Vai-5} {\sc J. V\"{a}is\"{a}l\"{a}}, The free quasiworld.
Freely quasiconformal and related maps in Banach spaces,
\textit{Quasiconformal geometry and dynamics $($Lublin 1996$)$,
55--118,
Banach Center Publ. {\bf 48}, Polish Acad. Sci. Inst. Math., Warsaw, 1999.}

\bibitem{Vai04}  {\sc J. V\"{a}is\"{a}l\"{a}}, Broken tubes in Hilbert spaces, \textit{Analysis (Munich)} {\bf 24} (2004), 227--238.

\bibitem{Vu1}  {\sc M. Vuorinen}, Capacity densities and angular limits of quasiregular mappings, \textit{Trans. Amer. Math. Soc.} {\bf 263} (1981), 343--354.

%\bibitem{Vu2}  {\sc M. Vuorinen,} Conformal invariants and quasiregular mappings, \textit{J. Anal. Math.} {\bf 45} (1985), 69--115.

%\bibitem{WZGHR}  {\sc X. Wang, Q. Zhou, T. Guan, M. Huang, and A. Rasila}, A note on the relationship between quasi-symmetric mappings and $\varphi$-uniform domains, \textit{J. Math. Anal. Appl.} {\bf 445} (2017), 1114--1119.


\bibitem{ZR} {\sc Q. Zhou and A. Rasila}, Quasim\"obius invariance of uniform domains, In press, \textit{Studia Math.} DOI: 10.4064/sm191215-22-10.
%\bibitem{ZRL} {\sc Q. Zhou, A. Rasila, and Y. Li}, Free distance ratio mappings in Banach spaces, \textit{Monatsh. Math.}  {\bf 191} (2020), 843--856.

\end{thebibliography}
\end{document}